\documentclass[a4paper]{article}
\usepackage{graphicx}
\usepackage{amssymb, amsmath, amsfonts, amsthm,amsbsy}
\usepackage[latin1]{inputenc}

\numberwithin{equation}{section}
\numberwithin{figure}{section}

\newcommand{\abs}[1]{\left\vert#1\right\vert}
\newcommand{\e}{\ensuremath{\mathrm{e}}}
\newcommand{\ii}{\ensuremath{\mathrm{i}}}
\newcommand{\dt}{\tau}
\newcommand{\dx}{h}
\newcommand{\ve}{\ensuremath{\varepsilon}}
\newcommand{\ol}[1]{\overline{#1}}

\title{Plane wave stability of some conservative schemes for the cubic Schr\"{o}dinger equation}
\author{Morten Dahlby and Brynjulf Owren}
\date{April 2, 2009}

\begin{document}
\maketitle

\begin{abstract}
	The plane wave stability properties of the conservative schemes of Besse and Fei et al.\ for the cubic Schr\"{o}dinger equation are analysed. Although the two methods possess many of the same conservation properties, we show that their stability behaviour is very different.
	An energy preserving generalisation of the  Fei method with improved stability is presented.
\end{abstract}
 
\section{Introduction }
	\begin{figure}[!ht]
		\centering
		\includegraphics[width=\textwidth]{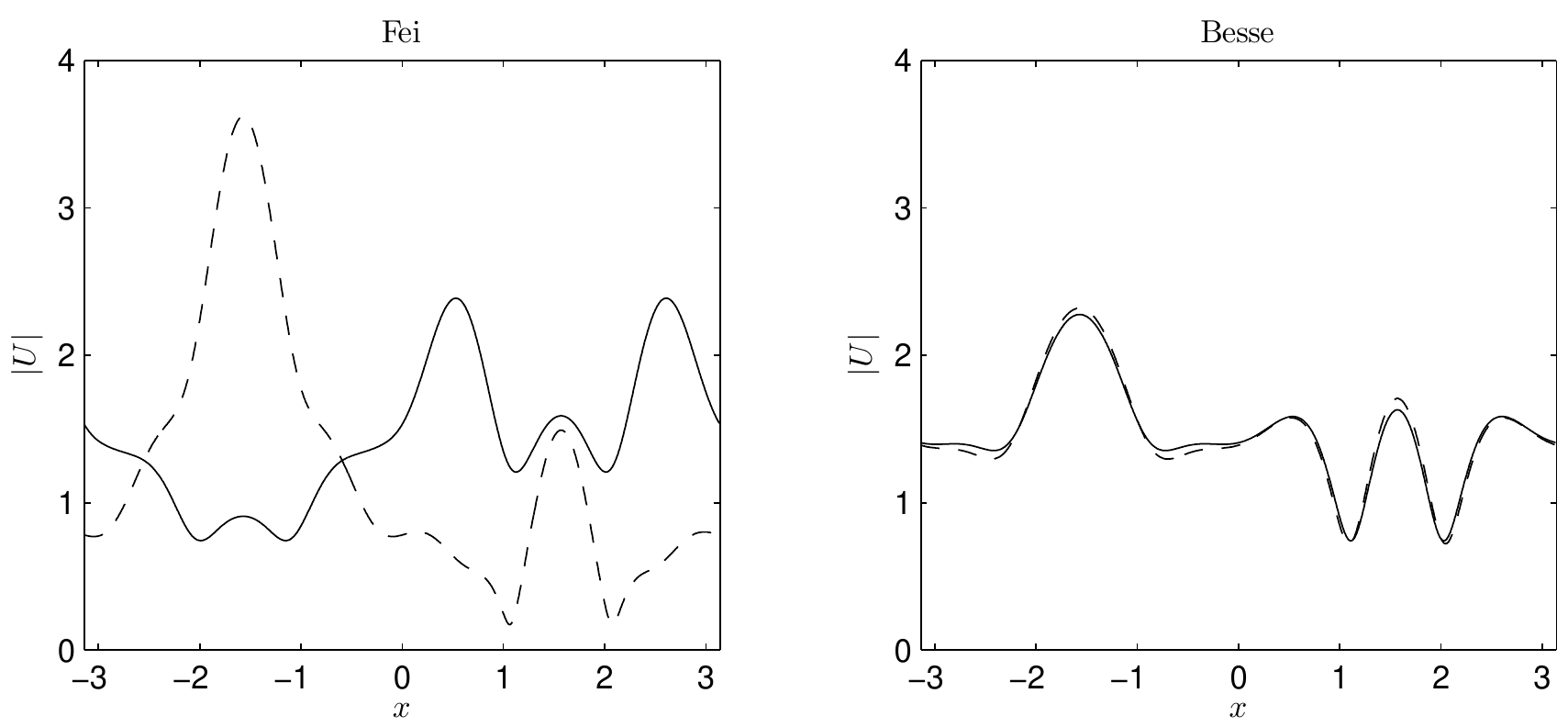}
		\caption{The plot shows the numerical solution of the cubic Schr\"{o}dinger equation produced by the method of Fei et al.\ \cite{fei95nso} and the method of Besse \cite{besse04ars}, in two consecutive time steps. Here $u_0(x)=\e^{\sin x},\ t=1.9$, time step $\dt=0.01$, space step $\dx=2\pi/1024,\ \lambda=2$.}
		\label{fig:feijump}
		\centering
	\end{figure}
	Figure~\ref{fig:feijump} shows the numerical approximation in two consecutive time steps with two well-known schemes  applied to the cubic  Schr\"{o}dinger equation (CSE).
	\begin{equation}\label{eq:nls}%
	   \ii u_t + \Delta u = \lambda|u|^2u,\quad\lambda\in\mathbb{R}.
	\end{equation}	
	Both schemes have order 2, are symmetric, and conserve mass and energy in a way to be made precise below. They are both linearly implicit. Whereas the method of Besse exhibits a continuous behaviour, the method of Fei et al.\ has an unstable spurious solution with eigenvalue near $-1$ which causes the numerical solution to alternate in each time step.

	We begin by briefly describing some properties of the CSE and some numerical methods which have been proposed for approximating its solution.
	In one space dimension, $d=1$, 
	\eqref{eq:nls} is completely integrable and for the sake of notational simplicity we shall assume in the rest of this paper that there is just one space dimension. All the numerical methods we consider can be easily generalised to arbitrary $d>1$. We may also add that the integrability which is particular to the case $d=1$ will not be important for the issues discussed here.
	Generally, the density
	\begin{equation} \label{eq:density}
	   \rho[u] = \int |u|^2\,\mathrm{d}x,
	\end{equation}
	is a conserved quantity, and it is also usual to work with the Hamiltonian structure
	\begin{equation} \label{eq:energy}
	    H[u] = \int \left(\frac{1}{2} |\nabla u|^2 + \frac{\lambda}{4} |u|^4\right)\,
	    \mathrm{d}x.
	\end{equation}
	In the literature, there exists a large variety of numerical schemes for the integration of \eqref{eq:nls}, see for instance 
	\cite{ablowitz76and, berland06stn,besse04ars, celledoni08sei,duran00tni,fei95nso,islas01gif,matsuo01doc,taha84aanII,weideman86ssm}. Frequently one sees examples of numerical schemes that are conservative, by this we mean that they exactly preserve some discretised version of
	\eqref{eq:density}, \eqref{eq:energy}, or both. In addition, it is often desirable to have schemes that are symmetric, see e.g.\ Hairer et al.\ \cite{hairer06gni} for a proper definition.
	One favourable consequence of having a conservative scheme, is that it can be used to control the growth of the numerical solution over long times. Generally, most conservative schemes are implicit in the time step. But the class of implicit schemes can be further divided into schemes that are fully implicit, and those that are linearly implicit or semi-explicit.
	In the former case, a system of nonlinear equations must be solved in every time step, the size of which is equal to the number of spatial degrees of freedom. It is often necessary to use a Newton-type iteration for solving this nonlinear system, and especially for large time steps the convergence may be slow or the iteration may not converge at all. However, using a linearly implicit scheme guarantees that the cost is roughly the same in each time step, therefore such schemes are attractive from the point of view of computational efficiency. Examples of linearly implicit schemes which are symmetric and conserve some discretised version of the energy when applied to the cubic Schr\"odinger equation are the method of Besse \cite{besse04ars} and that of Fei et al.\ \cite{fei95nso}, the latter being derived also in \cite{matsuo01doc}. The former scheme is formulated continuously in space, i.e.\ for each time $0=t_0<t_1<\cdots$ the method produces an approximation $U^n$ to $u(x,t_n)$ as follows
	\begin{align}
	\frac{\phi^{n+1/2}+\phi^{n-1/2}}{2} &= |U^n|^2, \label{eq:besse1} \\[2mm]
	\ii\,\frac{U^{n+1}-U^n}{\dt}+
	\Delta\left(
	\frac{U^{n+1}+U^n}{2}
	\right)&=\lambda \phi^{n+1/2}\,\left(
	  \frac{U^{n+1}+U^n}{2}
	\right). \label{eq:besse2}
	\end{align}
	This method can alternatively be written as a two-step scheme, and thus an auxiliary approximation is needed for the first step. Besse proposes to set
	$(U^0,\phi^{-1/2})=(u(x,0),|u(x,0)|^2)$, and shows that the following two versions of \eqref{eq:density} and \eqref{eq:energy}
	\begin{equation*}\label{eq:bessedensity}
	D^n = \int |U^n|^2\,\mathrm{d}x,
	\end{equation*}
	\begin{equation*} \label{eq:besseenergy}
	H^n = \int\left(\frac{1}{2}|\nabla U|^2 + \frac{\lambda}{4}\phi^{n+1/2}\phi^{n-1/2}\right)\, \mathrm{d}x,
	\end{equation*}
	are preserved in the sense that $D^n=D^0$ and $E^n=E^0$ for all $n$.

	The method of Fei et al.\ is different from the Besse method in that it is discretised both in time and space, in
	\cite{fei95nso} the scheme is given for one space dimension as follows:
	\begin{equation} \label{eq:fei}
	\ii \frac{U_m^{n+1}-U_m^{n-1}}{2\dt} +
	\frac{1}{\dx^2}\delta_{x}^2\left(\frac{U_m^{n+1}+U_m^{n-1}}2\right)=\lambda\left|U_m^n\right|^2\left(
	\frac{U_m^{n+1}+U_m^{n-1}}2\right),
	\end{equation}
	where $\delta_{x}$ is the centered difference in  space, that is,
	$\delta_{x} U_m^n=U_{m+\frac{1}{2}}^n-U_{m-\frac{1}{2}}^n$. Here $U_m^n$ is an approximation to $u(x_m,t_n)$. The associated discretised density and energy of this scheme are
	\begin{equation*}\label{eq:feidensity}
	D^n = \frac{h}{2}\sum_m \left( |U_m^n|^2 + |U_m^{n+1}|^2\right),
	\end{equation*}
	\begin{equation*}\label{eq:feienergy}
	H^n = \frac{1}{4}h\sum_m\left(  \left|\frac{U_{m+1}^{n+1}-U_m^{n+1}}{h}\right|^2
	+\left|\frac{U_{m+1}^{n}-U_m^{n}}{h}\right|^2
	+\lambda  |U_{m}^{n+1}|^2|U_{m}^{n}|^2\right).
	\end{equation*}
	In a similar way as for the Besse scheme, \eqref{eq:fei} can be formulated continuously in space as
	\begin{equation} \label{eq:feicontspace}
	 \ii \frac{U^{n+1}-U^{n-1}}{2\dt} +
	 \Delta\left(\frac{U^{n+1}+U^{n-1}}{2}\right)=\lambda\left|U^n\right|^2\left(
	\frac{U^{n+1}+U^{n-1}}2\right).
	\end{equation}

	The rest of the paper is organised as follows; in section~\ref{sec:analysis} we analyse the plane wave stability properties of the two schemes 
	\eqref{eq:besse2} and \eqref{eq:fei}. We illustrate the results using various plots.  
	In section~\ref{se:mod} we introduce a two parameter energy conserving generalisation of the Fei method that improves stability. 


\section{Plane wave solutions, dispersion relations and stability}
		\label{sec:analysis}
		Most of this section is inspired by the analysis of \cite{weideman86ssm} where dispersion relations for the exact and some
		numerical solutions of \eqref{eq:nls} are derived and their linear stability is analysed. 
		\subsection{The exact solution}
		The cubic Schr\"{o}dinger equation \eqref{eq:nls} with periodic boundary conditions supports plane wave solutions of the form
		\begin{equation*}
		     v(x,t) = a\e^{\ii (kx-\omega t)},
		\end{equation*}
		where $\omega$ is determined by the dispersion relation
		\begin{equation} \label{eq:disprelnls}
		  \omega = k^2 + \lambda |a|^2.
		\end{equation}
		We then consider small perturbations of such a solution, substituting
		$u(x,t)=(1+\ve(x,t))v(x,t)$ into \eqref{eq:nls}, and ignoring $\mathcal{O}(\ve^2)$ terms. We get
		\begin{equation} \label{deq:eps}
		    \ii\ve_t +2\ii k\ve_x +\ve_{xx} -\lambda |a|^2(\ve+\ol{\ve}) = 0.
		\end{equation}
		The perturbation $\ve(x,t)$ is periodic, thus we invoke its Fourier expansion
		\begin{equation*}
		  \ve(x,t) = \sum_{\ell\in\mathbb{Z}} \hat{\ve}_\ell(t)\e^{\ii\ell x}.
		\end{equation*}
		 Substitute this expansion into \eqref{deq:eps} to obtain
		 \begin{equation*}
		 \frac{\mathrm{d}}{\mathrm{d}t}\hat{\ve}_\ell =
		 \ii\left( (-2k\ell-\ell^2-\lambda |a|^2)\hat{\ve}_\ell-\lambda|a|^2\ol{\hat{\ve}}_{-\ell}
		 \right).
		 \end{equation*}
		We take the complex conjugate of this last equation and replace $\ell$ by $-\ell$, the result is the linear system of ODEs
		\begin{equation*}\label{eq:ODE2by2}
			\frac{\mathrm{d}}{\mathrm{d}t}
			\begin{bmatrix}
				\hat{\varepsilon}_{\ell}\\
				\ol{\hat{\varepsilon}}_{-\ell}			
			\end{bmatrix}
			=
			\ii G_\ell
			\begin{bmatrix}
				\hat{\varepsilon}_\ell\\
				\ol{\hat{\varepsilon}}_{-\ell}
			\end{bmatrix},\quad
			G_\ell=
			\begin{bmatrix}
				-\ell^2-2k\ell-\lambda|a|^2 & -\lambda|a|^2\\
				\lambda|a|^2 & 	\ell^2-2k\ell+\lambda|a|^2
			\end{bmatrix}.
		\end{equation*}
		The eigenvalues of $G_\ell$ are
		\begin{equation*}\label{eq:Geigval}%
		\lambda_\ell = \left(-2k\pm\sqrt{\ell^2+2\lambda|a|^2}
		    \right)\ell.
		\end{equation*}
		If $\lambda_\ell$ is complex, then one eigenvalue of $(\ii\,G_\ell)$
		will have a positive real part, and the corresponding mode will be
		unstable. This happens if
		\begin{equation*}
		\ell^2 < -2\lambda|a|^2,
		\end{equation*}
		thus instability may only occur when $\lambda < 0$, usually referred to as the focusing case.

	\subsection{The scheme of Fei et al.}\label{se:fei}
		We now consider again the scheme 
			\begin{equation} \label{eq:fei1d}
			\ii \frac{U_m^{n+1}-U_m^{n-1}}{2\dt} +
			\frac{1}{\dx^2}\delta_x^2\left(\frac{U_m^{n+1}+U_m^{n-1}}2\right)=\lambda\left|U_m^n\right|^2\left(
			\frac{U_m^{n+1}+U_m^{n-1}}2\right).
		\end{equation}
		Substituting a sequence of the form $V_m^n = a\e^{\ii(kx_m-\omega t_n)}$
		in which $x_m = m\dx$, $t_n = n\dt$, we get the dispersion relation
		\begin{equation} \label{eq:feidisp}
			\tan\omega\dt = \lambda\tau |a|^2 + 4\rho\sin^2\frac{kh}{2},\quad
			\rho=\frac{\dt}{\dx^2}.
		\end{equation}
		We perturb this plane wave solution, substituting $U_m^n=(1+\ve_m^n)V_m^n$ into \eqref{eq:fei1d} and after ignoring quadratic and higher order terms in $\ve_m^n$ we get
		\begin{equation} \label{eq:feipertdiff}
			a_1\varepsilon_{m+1}^{n+1}+a_0\varepsilon_m^{n+1}+a_{-1}\varepsilon_{m-1}^{n+1}
			= b(\varepsilon_m^n + \ol{\varepsilon}_m^n) -
			\ol{a}_{-1}\varepsilon_{m+1}^{n-1}-\ol{a}_0\varepsilon_m^{n-1}-\ol{a}_1\varepsilon_{m-1}^{n-1} .
		\end{equation}		
		Here
		\begin{align*}
			a_1 &= \rho\, \e^{\ii(k\dx-\omega \dt)}& b&=2q\cos\omega\dt\\
			a_0 &= \e^{-\ii\omega\dt}\,(\ii-2\rho-q)& q&=\lambda\dt|a|^2\\
			a_{-1} &=  \rho\, \e^{-\ii(k\dx+\omega \dt)}.&&
		\end{align*}
		We expand $\ve_m^n$ in a Fourier series 
		$\ve_m^n =\displaystyle{ \sum_{\ell\in\mathbb{Z}}} \hat{\ve}_\ell^n\e^{\ii\ell x_m}$
		and substitute this into \eqref{eq:feipertdiff} to obtain
		\begin{equation*}
		c_\ell\hat{\ve}_\ell^{n+1}=b(\hat{\ve}_\ell^n+\ol{\hat{\ve}}_{-\ell}^n)
		-\ol{c}_\ell\hat{\ve}_\ell^{n-1},\qquad
		c_\ell = a_1\e^{\ii \ell\dx}+a_0+a_{-1}\,\e^{-\ii \ell\dx}.
		\end{equation*}
		We now take the complex conjugate of this equation and replace $\ell$ by $-\ell$
		to obtain a system of difference equations for $E^n=\left[\hat{\ve}_\ell^n,\ol{\hat{\ve}}_{-\ell}^n\right]^T$,
		\begin{equation*}
			\begin{bmatrix}
				c_\ell & 0 \\
				0   & \ol{c}_{-\ell}
			\end{bmatrix}E^{n+1}=
			\begin{bmatrix}
				b & b \\
				b & b			
			\end{bmatrix}E^n-
			\begin{bmatrix}
				\ol{c}_\ell & 0 \\
				0   & c_{-\ell}			
			\end{bmatrix}E^{n-1}.
		\end{equation*}
		We find that this difference equation is stable if and only if
		the polynomial
		\begin{equation} \label{eq:feistabpol}
			p(z)=
			c_\ell\ol{c}_{-\ell}z^4
			-b(c_\ell + \ol{c}_{-\ell}) z^3
			+(c_\ell c_{-\ell}+\ol{c}_\ell\ol{c}_{-\ell})z^2
			-b(c_{-\ell}+\ol{c}_{\ell})z
			+c_{-\ell}\ol{c}_\ell
		\end{equation}
		has all its roots in the closed unit disc. Note that $p(z)$ is \emph{self-reciprocal}, 
		meaning that its set of roots is invariant under the transformation $z\mapsto 1/\ol{z}$.
		Each root on the unit circle is invariant under this transformation, but any root in the open unit disc is mapped to a root outside the unit circle. Thus, for self-reciprocal polynomials, stability is equivalent to all roots lying on the unit circle.
		 Note that when we use the scheme which is continuous in space \eqref{eq:feicontspace}, we
		 get again the stability polynomial \eqref{eq:feistabpol}, but where 
		 $c_\ell$ is replaced by
		 \begin{equation*}
		c_\ell = \e^{-\ii\omega \dt}\left(\ii-(K+L)^2-q\right),\qquad
		K=k\sqrt{\tau},\ L=\ell\sqrt{\tau}.
		\end{equation*}
		The case $k=0$ is particularly simple, since the stability polynomial in this case will have real coefficients. In that case one may easily derive that $p(z)$ has all its root on the unit circle for all values of $\ell$ if and only if 
		\begin{equation*}
		\lambda |a|^2 \leq \frac{1-\cos\dx}{\dx^2}.
		\end{equation*}
		Note that this condition is independent of $\dt$, and in the limit when $\dx$ tends to zero, one simply gets the condition $\lambda |a|^2\leq \frac{1}{2}$. This does however not imply that the scheme does not converge on finite time intervals $[0,t^*]$. Suppose that $\lambda |a|^2 > \frac{1}{2}$.
		For small perturbations, one may expect that the error is roughly amplified with a factor in each step that equals the magnitude of the largest root of \eqref{eq:feistabpol}, which for $k=0$ and in the limit case $\dx\rightarrow 0$ is of the form
		\begin{equation*}
		z_* = -\left(1+\tau \ell\sqrt{2\lambda |a|^2-\ell^2}\right) + \mathcal{O}(\tau^2),
		\end{equation*}
		thus locally the size of the spurious solution grows exponentially, the growth factor of the unit frequency over the interval $[0,t^*]$ being approximately $\exp(C t^*)$ with $C=\sqrt{2\lambda |a|^2-1}$.

	\subsection{The scheme of Besse}
		In one space dimension the scheme of Besse \eqref{eq:besse1}, \eqref{eq:besse2} is
		\begin{align} \label{eq:besse1da}
			\frac{\phi^{n+\frac12}+\phi^{n-\frac12}}2&=\abs{U^n}^2, \\ \label{eq:besse1db}
			\frac{\ii}{\dt}\left( U^{n+1}-U^{n}\right) +
			\frac{U^{n+1}_{xx}+U^{n}_{xx}}2&=\lambda\phi^{n+\frac12}\left(\frac{U^{n+1}+U^{n}}2\right).
		\end{align}
		We use the plane wave solution $V^n=a\e^{\ii(kx-\omega t_n)}$, which now is continuous in space, to get the following relation
		for $\omega$.
		\begin{equation}\label{eq:bessedisp}
			\tan{\frac{\omega\dt}2}=\frac12\left(\lambda\dt\abs{a}^2+\dt k^2\right).
		\end{equation}
		As for the Fei scheme we consider the perturbations
		\begin{equation*}
			U^n=V^n\left(1+\ve^n\right),\quad\phi^{n+\frac12}=\abs{a}^2\left(1+\delta^{n+\frac12}\right).
		\end{equation*}
		Substituting these expressions into \eqref{eq:besse1da} and ignoring higher order terms yields
		\begin{equation}\label{eq:bessediffa}
			\delta^{n+\frac12}=-\delta^{n-\frac12}+2\left(\ve^n+\ol{\ve}^n\right).
		\end{equation}
		We now plug the last three expressions into
		  \eqref{eq:besse1db} to obtain
		\begin{multline}\label{eq:bessediffb}
			\frac{\ii}{\dt}\left(\e^{-\ii\omega\tau}\ve^{n+1}-\ve^n\right)+\frac12 \e^{-\ii\omega\tau}\left(\ve^{n+1}_{xx}+2\ii k\ve^{n+1}_x-\left(k^2+\lambda\abs{a}^2\right)\ve^{n+1}\right)\\
			+\frac{1}2\left(\ve^{n}_{xx}+2\ii k\ve^{n}_x-\left(k^2+\lambda\abs{a}^2\right)\ve^{n}\right)=\frac12\lambda\abs{a}^2\left(\e^{-\ii\omega\tau}+1\right)\delta^{n+\frac12}.
		\end{multline}
		Expand $\ve^n(x)$ and $\delta^{n+\frac12}(x)$ in a Fourier series and insert
		into \eqref{eq:bessediffa} and \eqref{eq:bessediffb} to get the system
		\begin{equation*}
			\begin{bmatrix}
				c_\ell&0&-b&0\\
				0&\ol{c}_{-\ell}&0&-b\\
				0&0&1&0\\
				0&0&0&1
			\end{bmatrix}E^{n+1}
			=
			\begin{bmatrix}
				-\ol{c}_\ell&0&0&0\\
				0&-c_{-\ell}&0&0\\
				2&2&-1&0\\
				2&2&0&-1
			\end{bmatrix}E^{n},
			\quad E^n=
			\begin{bmatrix}
				\hat{\ve}^n_{\ell}\\
				\ol{\hat{\ve}}^n_{\ell}\\
				\hat{\delta}_{\ell}^{n-\frac12}\\
				\ol{\hat{\delta}}_{\ell}^{n-\frac12}
			\end{bmatrix},
		\end{equation*}
		where we have defined
		\begin{align*}
			c_\ell&=(2\ii-(L+K)^2-q)\e^{-\ii\omega\tau}\\
			b&=2q\cos\frac{\omega\dt}2,
		\end{align*}
		and where $L=\sqrt{\tau}\ell$ and $K=\sqrt{\tau}k$ as before. 
		Notice that $c_\ell$ and $b$ are defined differently than in the Fei case.
		We find that the characteristic polynomial is $(z+1)\tilde{p}(z)$ with
		\begin{multline*}
			\tilde{p}(z)=c_\ell\ol{c}_{-\ell}z^3+\left(c_\ell c_{-\ell}+\ol{c}_\ell c_{-\ell}+c_\ell\ol{c}_{-\ell}-2b\left(c_\ell+\ol{c}_{-\ell}\right)\right)z^2\\
			+\left(c_\ell c_{-\ell}+\ol{c}_\ell \ol{c}_{-\ell}+\ol{c}_\ell c_{-\ell}-2b\left(c_\ell+\ol{c}_{-\ell}\right)\right)z+\ol{c}_\ell c_{-\ell}.
		\end{multline*}
		To shorten the notation we divide by $c_\ell\ol{c}_{-\ell}$ and define $f$ and $g$ such that 
		\begin{equation*}
			p(z)=z^3+\ol{g}z^2+gz+f,
		\end{equation*}
		and $p(z)$ has the same roots as $\tilde{p}(z)$. 
		The polynomial is self-reciprocal and we have stability only if all three roots lie on the unit circle. We proceed to 
		express the stability region $S$ in terms of $g$ for a given $f$. The key observation is that for points on the boundary $\partial S$ at least  
		 two of the roots are coalescing. For such values of $f$ and $g$, we can write the polynomial as
		\begin{equation*}
				p(z)=\left(z-\e^{\ii\theta}\right)^2\left(z-\e^{\ii\psi}\right).
		\end{equation*}
		By comparing coefficients we get that for a given value of $f$ we can parametrise the stability boundary for $g$ as follows
		\begin{equation*}
			g(\theta)=\e^{2\ii\theta}-2f\e^{-\ii\theta}.
		\end{equation*}
		Since $\abs{f}=1$ we have that $\abs{g}\leq 3$ is necessary for stability while $\abs{g}\leq 1$ is sufficient. 
		For $k=0$ the latter condition becomes $\ell^2<-2\lambda|a|^2$, which is the same as for the exact solution.
		
		For $k=0$ the two roots of $p(z)$ with the largest magnitude are on the form
		\begin{equation*}
			z_*=1\pm\tau\ell\sqrt{-2\lambda|a|^2-\ell^2}+\mathcal{O}(\tau^2).
		\end{equation*}
		In this case the eigenvalue is near $1$, not $-1$ as in the Fei case, which explains why the Besse scheme does not 
		exhibit an alternating behaviour for $k=0$.

\paragraph{Extension to more space dimensions.}
The scheme of Fei \eqref{eq:fei1d} has an obvious generalisation to $d$ space dimensions with corresponding energy and density functions. We may also consider the general form of the Besse scheme \eqref{eq:besse1}, \eqref{eq:besse2}, and introduce plane wave solutions in $d$ dimensions
\[
u(\boldsymbol{x},t) = a\e^{\ii\left(\boldsymbol{k}\cdot\boldsymbol{x}-\omega t\right)}
\]
for space variables $\boldsymbol{x}=(x_1,\ldots,x_d)$ and wave numbers
$\boldsymbol{k}=(k_1,\ldots,k_d)$. Exact solutions of this form satisfy the dispersion relation
\[
       \omega = |\boldsymbol{k}|^2+\lambda |a|^2,
\]
and these solutions are stable with respect to perturbations
\[
   \ve_{\boldsymbol\ell}(\boldsymbol{x},t)=
   \hat{\ve}_{\boldsymbol\ell}(t)
   \e^{\ii\,\boldsymbol \ell\cdot\boldsymbol x},
\]
whenever $|\boldsymbol\ell|^2\geq -2\lambda |a|^2$.

For the method of Fei et al.\ one finds again that for $\boldsymbol k = \mathbf{0}$ the scheme is stable to perturbations only if $\lambda |a|^2 \leq \frac{1}{2}$ in the limit $h\rightarrow 0$, the critical perturbation wavenumber vectors being the canonical unit vectors in $\mathbb{R}^d$. The dispersion relation is now
\[
\tan\omega\tau = \lambda\tau|a|^2+4\rho\sum_{j=1}^d \sin^2\frac{k_j h}{2}\quad
\stackrel{h=0}{\longrightarrow}\quad
\lambda\tau |a|^2+\tau|\boldsymbol{k}|^2.
\]
Also the method of Besse generalizes similarly in more space dimensions, the dispersion relation being
\[
\tan\frac{\omega\tau}{2}=\frac{1}{2}(\lambda\tau |a|^2+\tau|\boldsymbol{k}|^2),
\]
and the stability polynomial is obtained simply by replacing $c_\ell$ by
\[
 c_{\boldsymbol\ell} = (2\ii - \tau|\boldsymbol \ell+\boldsymbol k|^2-dq)\e^{-\ii\omega\tau}.
\]

	\subsection{Numerical results}\label{se:plots}
\begin{figure}
			\centering
			\includegraphics[width=0.6\textwidth]{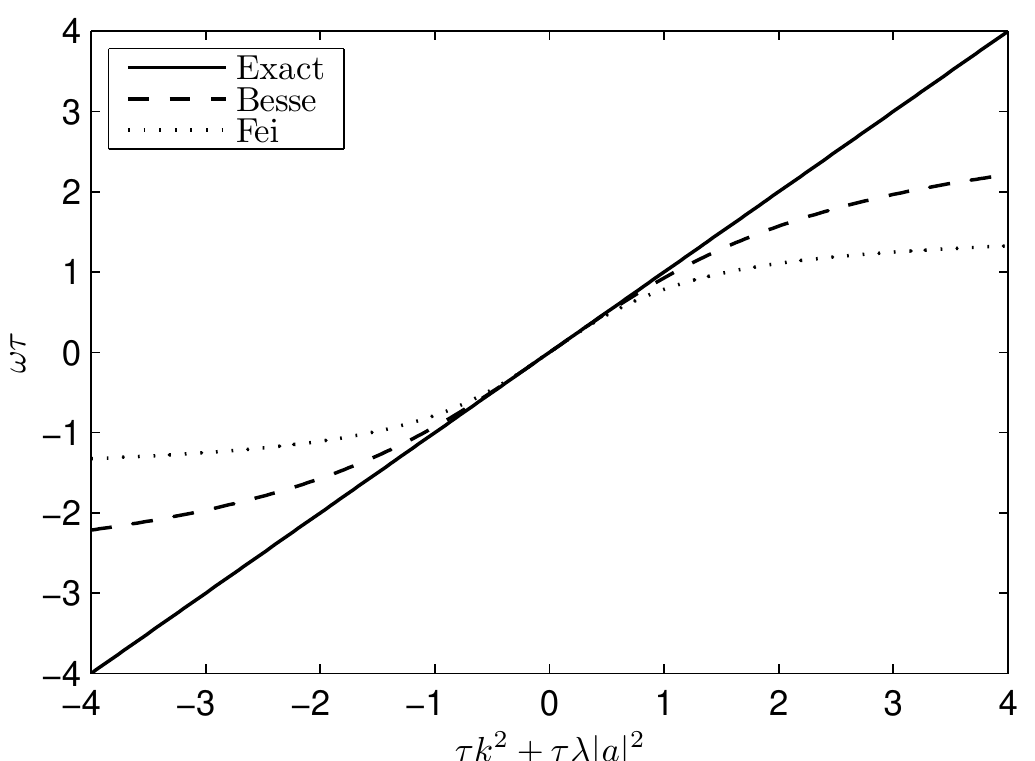}
			\caption{A comparison of the dispersion relations.}
			\label{fi:disprel}
		\end{figure}

 We compare the schemes of Besse and Fei et al.\ in the limit $\dx\rightarrow 0$, in which case $p(z)$ depends on the three parameters $K, L$ and $q$.
		In figure~\ref{fi:disprel} we compare the exact dispersion relation \eqref{eq:disprelnls} with Fei \eqref{eq:feidisp}
		and Besse \eqref{eq:bessedisp}. 
	         Both relations express $\omega\tau$ in terms of
	         $k^2+\lambda |a|^2$ and for the Besse method it is obtained by replacing the time step $\tau$ by $\tau/2$ in the Fei method.	
 
\begin{figure}
			\centering
			\includegraphics[width=\textwidth]{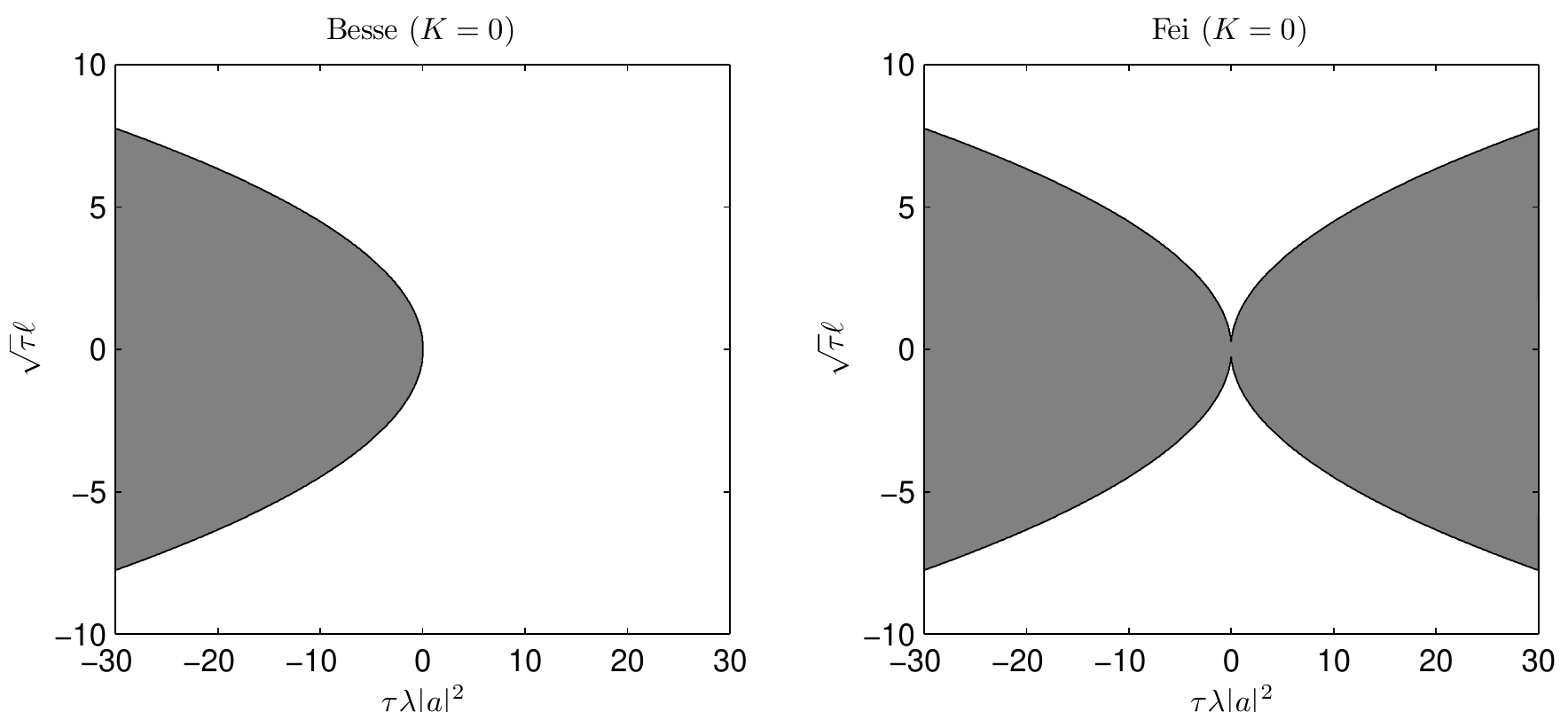}
			\caption{Stability for $K=0$ (grey is unstable).}
			\label{fi:stabplot_K0}
		\end{figure}

\begin{figure}
			\centering
			\includegraphics[width=\textwidth]{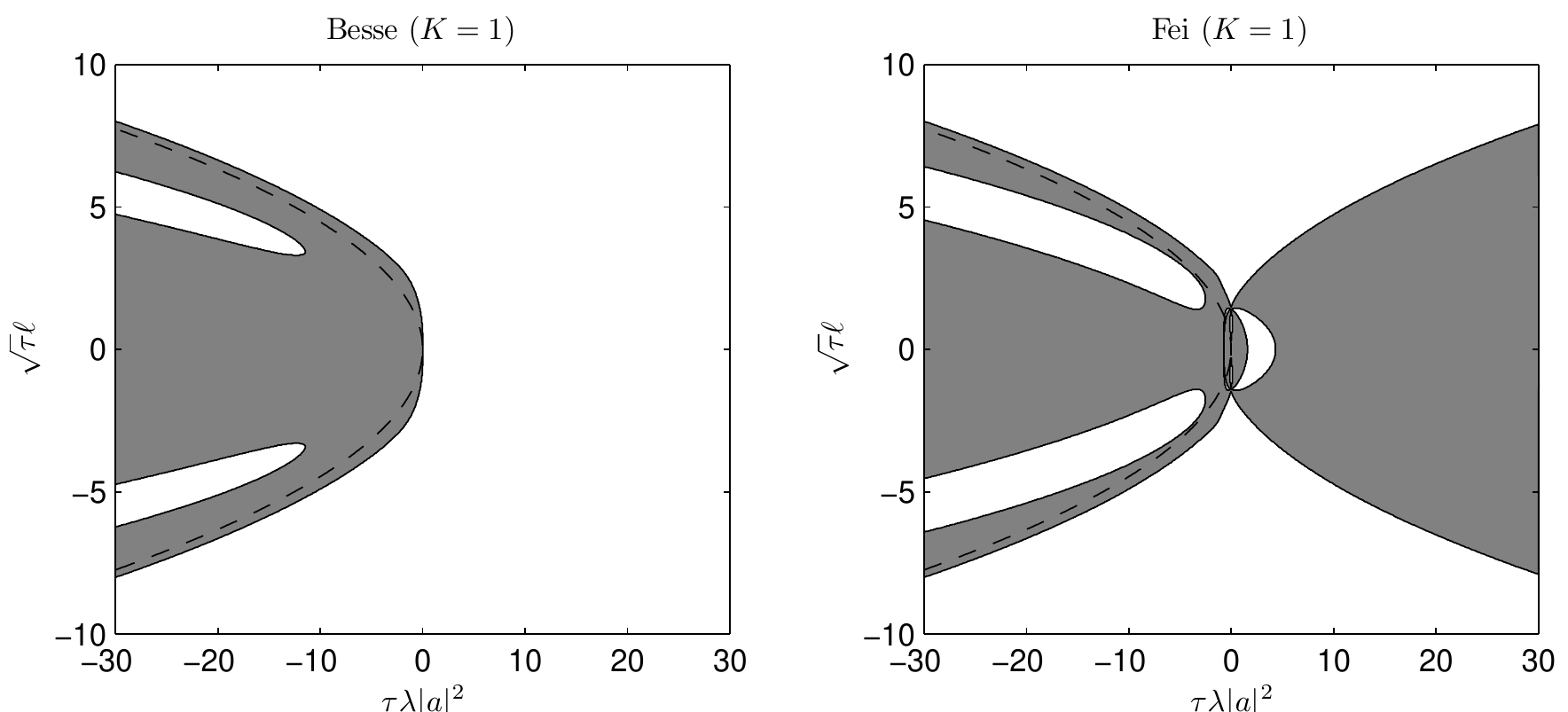}
			\caption{Stability for $K=1$.}
			\label{fi:stabplot_K1}
		\end{figure}

		Figure~\ref{fi:stabplot_K0} shows a numerical computation of the stability region for both Besse and Fei when fixing $K=k=0$ in $p(z)$. These plots are obtained from the explicit expressions in the analysis of the previous sections. The case $K=k\sqrt{\tau}=1$ is shown in figure~\ref{fi:stabplot_K1}. For reference, we have included the stability boundary of the exact solution as a broken line.
The stability of the Besse scheme differs slightly from that of the exact solution when $\lambda$ is negative, however, the stability is retained for positive $\lambda$. 
The Fei scheme is again unstable in the defocusing case. We have observed that when $K$ is further increased, then also the Besse method becomes unstable for small modes in the defocusing regime.

\begin{figure}
			\centering
			\includegraphics[width=0.6\textwidth]{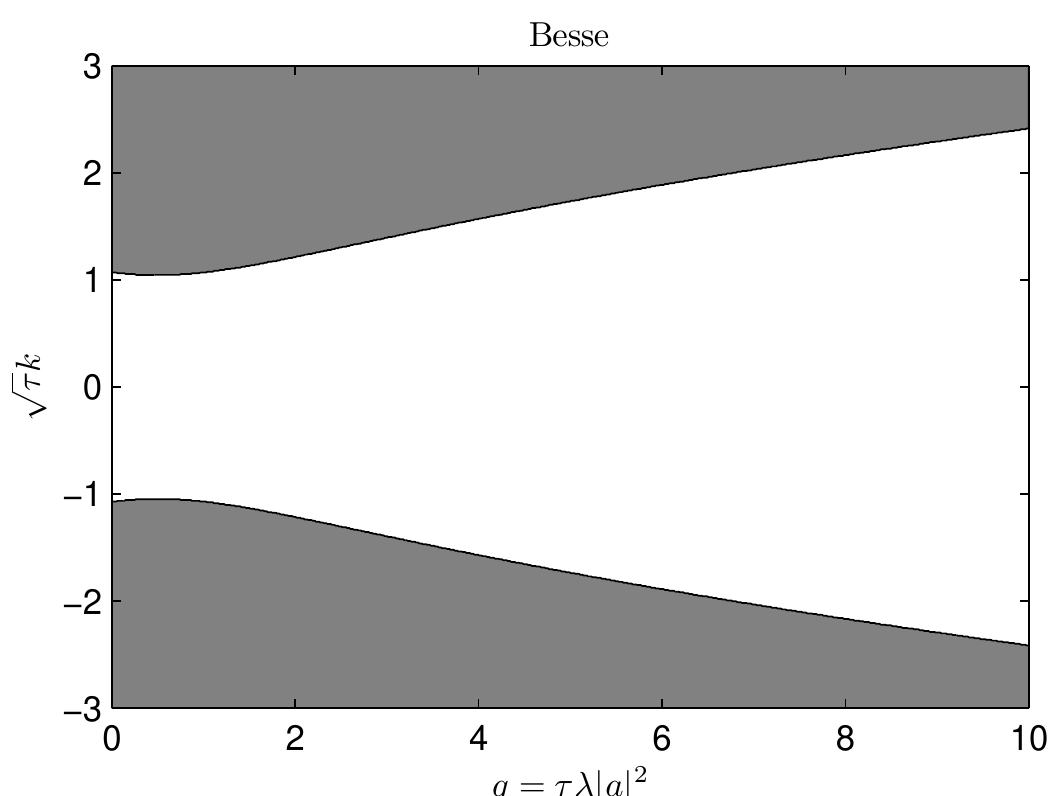}
			\caption{Stability for all $\ell$.}
			\label{fi:stabplot_besse}
		\end{figure}

		For a numerical plane wave solution to be stable, none of its  Fourier modes $\ell$ should be amplified by the method. The plot in figure~\ref{fi:stabplot_besse}
		shows the stability region in the $qK$-plane for the Besse method. A pair $(q,K)$ is unstable if
		the largest root of the stability polynomial exceeds 1 in modulus for some real    value of $L$. It is then only of interest to consider the defocusing case ($q\geq 0$) since the exact solution itself is unstable in this sense for all negative $q$. The Fei scheme is unstable for all $(q,K)$ where $q\neq 0$.
		
\section{An energy preserving modification of the Fei scheme}\label{se:mod}
	Using the procedure in \cite{matsuo01doc} one can derive the following symmetric 2-step scheme. 
	\begin{multline*}
		\ii \frac{U^{n+1}-U^{n-1}}{2\dt} +
		\left(\frac{\theta}{2}U_{xx}^{n+1}+
		(1-\theta)U_{xx}^n+\frac{\theta}{2}U_{xx}^{n-1}\right)\\
		=\frac{\lambda\gamma}2\left|U^n\right|^2\left(
		U^{n+1}+ U^{n-1}\right)+\frac{\lambda(1-\gamma)}2 \left(U^n\right)^2\left(\ol{U}^{n+1}+\ol{U}^{n-1}\right).
	\end{multline*}
	$\theta$ and $\gamma$ are real parameters. Note that $\theta=\gamma=1$ yields the method of Fei et al. This method is by construction 
	energy preserving, and its energy function is given as	
	\begin{align*}
	H^n&=\frac{\theta h}{4}\sum_m \left(|\delta^+U_m^{n+1}|^2+|\delta^+U_m^{n}|^2\right)\\
		&+\frac{(1-\theta)h}{4}\sum_m \left((\delta^+U_m^{n})(\delta^+\ol{U}_m^{n+1})+(\delta^+\ol{U}_m^{n})(\delta^+U_m^{n})\right)\\
		&+\frac{\gamma h}{4}\sum_m \lambda|U_m^{n}|^2|U_m^{n+1}|^2\\
		&+\frac{(1-\gamma)h}{16}\sum_m\lambda\left((U_m^{n+1})^2+(\ol{U}_m^{n+1})^2\right)\left((U_m^{n})^2+(\ol{U}_m^{n})^2\right)\\		
		&-\frac{(1-\gamma)h}{16}\sum_m\lambda\left((U_m^{n+1})^2-(\ol{U}_m^{n+1})^2\right)\left((U_m^{n})^2-(\ol{U}_m^{n})^2\right).		
	\end{align*}
	The last two lines comes from the relation $|u|^4=\frac14\left(u^2+\ol{u}^2\right)-\frac14\left(u^2-\ol{u}^2\right)$.
	We are only aware of a conserved density function in the Fei case, we will however see that some of the parameters yield
	improved stability compared to Fei. 
	
\begin{figure}
		\centering
		\includegraphics[width=\textwidth]{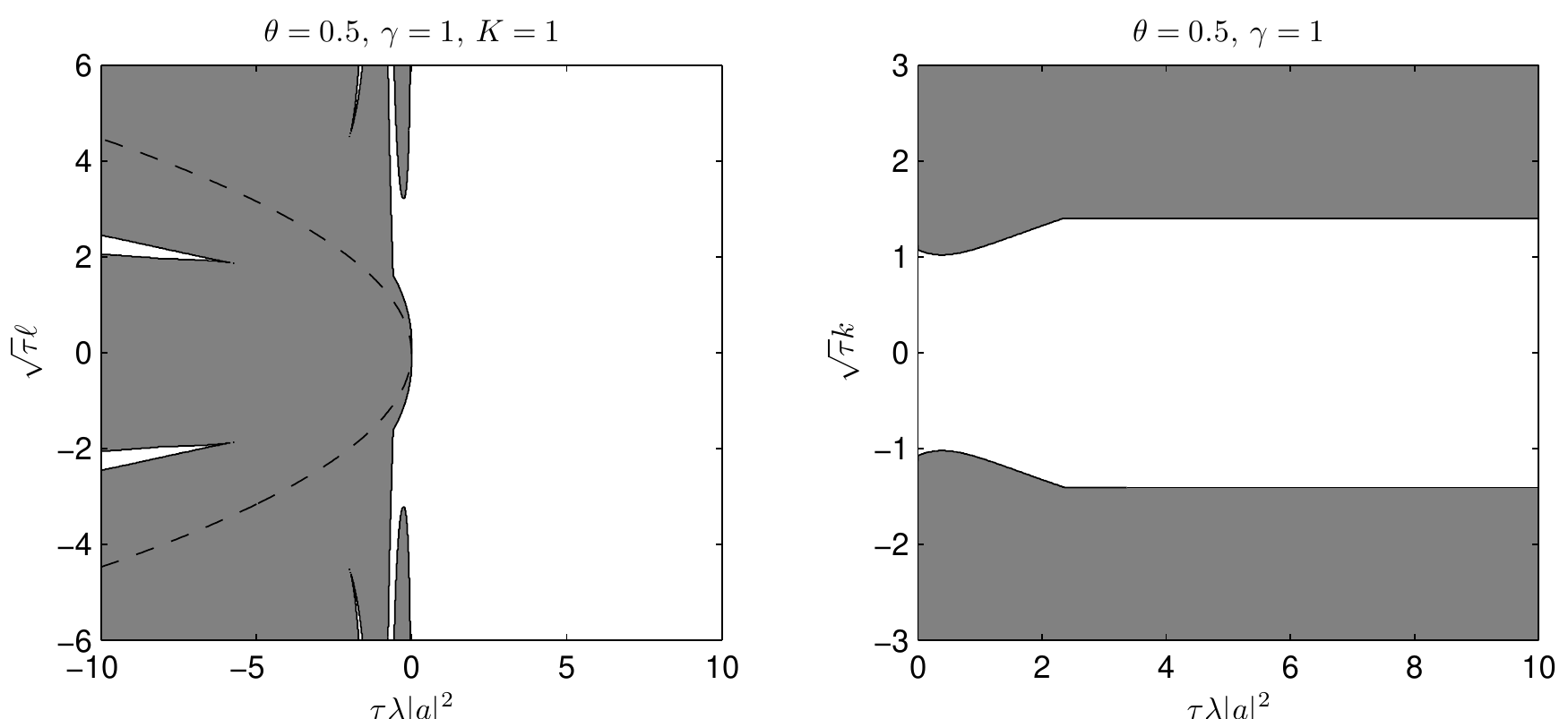}
		\caption{Stability for the modified scheme.}
		\label{fi:new}
	\end{figure}

	The dispersion relation is
	\begin{equation*}
		\sin\omega\tau=K^2\left(\theta\cos\omega\tau +(1-\theta)\right)+q\cos\omega\tau.
	\end{equation*}
	Using the same procedure as in chapter~\ref{se:fei} we get the stability polynomial
	\begin{align*}
		p(z)&=\left(c_ \ell\ol{c}_{- \ell}-(1-\gamma)^2q^2\right)z^4\\
		&+\left(\gamma(\gamma-1)b^2-\left((2-\gamma)b+d_ \ell\right)\ol{c}_{- \ell}-\left((2-\gamma)b+d_{- \ell}\right)c_ \ell\right)z^3\\
		&+\left(b(2-\gamma)(d_ \ell+d_{- \ell})+d_ \ell d_{- \ell}+c_ \ell c_{- \ell}+\ol{c}_ \ell\ol{c}_{- \ell}+2(1-\gamma)^2q(b-q)+4(1-\gamma)b^2\right)z^2\\
		&+\left(\gamma(\gamma-1)b^2-\left((2-\gamma)b+d_ \ell\right)c_{- \ell}-\left((2-\gamma)b+d_{- \ell}\right)\ol{c}_ \ell\right)z\\
		&+\left(\ol{c}_ \ell c_{- \ell}-(1-\gamma)^2q^2\right),
	\end{align*}
	where
	\begin{equation*}
		c_\ell= \left(\ii-\gamma q -\theta(K+L)^2\right)\e^{-\ii\omega\tau}\quad\text{and}\quad d_\ell=2(1-\theta)(K+L)^2.
	\end{equation*}
	The scheme is stable provided that all zeroes of the self-reciprocal polynomial $p(z)$ is on the unit circle. In figure~\ref{fi:new}
	we see that the stability region is almost the same size as for Besse, compare with figure~\ref{fi:stabplot_besse}.
	These parameters clearly yield an improvement over the Fei method. 

\section{Conclusion}
We have seen that two schemes which have essentially the same conserving properties and are both symmetric and linearly implicit show a 
rather different behaviour when applied to the cubic Schr\"odinger equation. This behaviour is analysed in terms of plane wave solutions.
Both schemes can be interpreted as two-step schemes and thus have a spurious solution, in the Fei scheme this solution is unstable even 
for small time steps contrary to the scheme of Besse. We show however that the scheme by Fei et al.\ can be stabilised without losing the 
symmetry or the energy preservation property.

It is interesting to observe that two different schemes, designed to be conservative in a very similar way, turns out behave completely differently with respect to stability. It remains to be seen if this phenomenon appears also in other PDE models than the cubic Schr\"{o}dinger equation.

 \bibliographystyle{plain}
 \bibliography{nlsbib}
 
 \end{document}